\documentclass[11pt]{amsart}
\usepackage[utf8]{inputenc}
\usepackage[T1]{fontenc}
\usepackage{lmodern}
\usepackage{amsmath,amssymb,amsthm}
\usepackage{graphicx}
\usepackage[a4paper,margin=2.8cm]{geometry}
\usepackage{enumitem}
\usepackage[hidelinks]{hyperref}

\newtheorem{theorem}{Theorem}[section]
\newtheorem{proposition}[theorem]{Proposition}
\newtheorem{lemma}[theorem]{Lemma}
\newtheorem{corollary}[theorem]{Corollary}
\theoremstyle{definition}
\newtheorem{definition}[theorem]{Definition}
\newtheorem{remark}[theorem]{Remark}

\newcommand{\HH}{\mathbb{H}}
\newcommand{\R}{\mathbb{R}}
\newcommand{\Id}{\mathrm{Id}}
\newcommand{\Inj}{\operatorname{Inj}}
\newcommand{\Argsh}{\mathit{Argsh}}
\renewcommand{\Im}{\operatorname{Im}}
\renewcommand{\Re}{\operatorname{Re}}
\newcommand{\cl}[1]{\overline{#1}}

\title[Irregular horocyclic orbits of infinite thinness]{On the closure of irregular orbits of the horocyclic flow with infinite thinness}
\author{Masseye Gaye and Amadou Sy}

\begin{document}

\begin{abstract}
The topological dynamics of the horocyclic flow $h_\R$ on the unit tangent bundle of a geometrically finite hyperbolic surface is well known: either $h_\R$ is minimal, or the minimal sets are the periodic orbits. When the surface is geometrically infinite, the situation is more complex and the possible presence of irregular orbits makes the description of minimal sets difficult. We construct a family of geometrically infinite hyperbolic surfaces $\Sigma_\delta$, depending on a real number $\delta$ and having an irregular horocyclic orbit of infinite asymptotic thinness whose closure is not minimal.
\end{abstract}

\maketitle


\section{Introduction}

In the study of a dynamical system, an important question is to understand its minimal sets, that is, the non-empty closed invariant sets which contain no proper non-empty closed invariant subset. For the horocyclic flow $h_\R$ on the unit tangent bundle of a geometrically infinite hyperbolic surface, restricted to its non-wandering set $\Omega_h$, this question is still open. Understanding the $h_\R$-minimal sets amounts to understanding the closure of irregular orbits (orbits that are neither closed nor dense in $\Omega_h$), see \cite{Bellis-these, Bellis, GayeLo, Kulikov, Matsumoto}.

For $u\in\Omega_h$, a classical argument to show that $\cl{h_\R(u)}$ is not minimal consists in establishing that
\[
T_u=\{t\in\R^*;\ g_t(u)\in\cl{h_\R(u)}\}\neq\emptyset,
\]
that is, the geodesic orbit $g_\R(u)$ meets the closure of $h_\R(u)$ at a non-zero time. This is related to the lower limit of the injectivity radius along the half-geodesic $u(\R_+)$, called the asymptotic thinness of $u(\R_+)$ and denoted $\Inj(u(\R_+))$. A.~Bellis showed \cite{Bellis}: \emph{if $\Sigma$ is geometrically infinite, $u\in\Omega_h$, $h_\R(u)$ is non-periodic and $0\le\Inj(u(\R_+))<+\infty$, then $\cl{h_\R(u)}$ is not $h_\R$-minimal; moreover $T_u=\R^*$ if $\Inj(u(\R_+))=0$, and $T_u$ contains an unbounded sequence if $0<\Inj(u(\R_+))<+\infty$.}

If $h_\R(u)$ is closed, then $T_u=\emptyset$. In \cite{Bellis-these}, Bellis asks whether $\Inj(u(\R_+))=+\infty$ always implies $T_u=\emptyset$, in other words whether infinite thinness gives rise to new $h_\R$-minimal sets. We show by an example that the answer is negative.

\begin{theorem}\label{thm:main}
Let $\delta$ be a real number satisfying $1<\delta<3$.
There exist a geometrically infinite hyperbolic surface $\Sigma_\delta$ and $u_0\in\Omega_h\subset T^1\Sigma_\delta$ such that:
\begin{enumerate}[label=(\arabic*)]
\item $h_\R(u_0)$ is an irregular orbit, $\Inj(u_0(\R_+))=+\infty$, and $\cl{h_\R(u_0)}$ is not $h_\R$-minimal;
\item $T_{u_0}$ contains the sequence $t_k=2k\ln\delta$, $k\ge1$, which tends to $+\infty$.
\end{enumerate}
\end{theorem}

Section~2 recalls the classification of limit points, the topological nature of horocyclic orbits and horocyclic convergence. Section~3 proves Theorem~\ref{thm:main}.

\section{Preliminaries}

\subsection{Nature of limit points.}
Let $\Gamma$ be a Fuchsian group without elliptic elements and $\Sigma=\Gamma\backslash\HH$. Let $\Lambda=\cl{\Gamma.i}\setminus\Gamma.i$ denote its limit set. Following \cite{DalboStarkov}, a point $\xi\in\Lambda$ is said to be
\begin{itemize}
\item \emph{horocyclic} ($\xi\in\Lambda_h$) if every open horodisc $O_\xi$ based at $\xi$ meets $\Gamma.i$ in an infinite set;
\item \emph{discrete} ($\xi\in\Lambda_d$) if every open horodisc $O_\xi$ meets $\Gamma.i$ in a finite set;
\item \emph{parabolic} ($\xi\in\Lambda_p$) if $\xi$ is fixed by a parabolic element of $\Gamma$;
\item \emph{irregular} ($\xi\in\Lambda_{irr}$) otherwise.
\end{itemize}
These four $\Gamma$-invariant sets form a partition of $\Lambda$.

\begin{remark}\label{rem:horo}
Suppose $\infty\in\Lambda$. The horodiscs based at $\infty$ are the sets $\{\Im z>s\}$, $s>0$. Consequently:
\begin{itemize}
\item $\infty\notin\Lambda_h$ if and only if there exists $M>0$ such that $\Im\gamma(i)\le M$ for all $\gamma\in\Gamma$;
\item if there exists a sequence of pairwise distinct elements $\gamma_n\in\Gamma$ with $\Im\gamma_n(i)\to l>0$, then $\infty\notin\Lambda_d$ (the horodisc $\{\Im z>l/2\}$ contains infinitely many points of the orbit).
\end{itemize}
\end{remark}

Recall \cite{Dalbo} that $\Sigma$ is geometrically finite if and only if $\Lambda=\Lambda_h\cup\Lambda_p$.

\subsection{Topological nature of horocyclic orbits.}
We consider $h_\R$ on $T^1\Sigma$ restricted to the non-wandering set $\Omega_h=\{u:\ \tilde u(\infty)\in\Lambda\}$. When $\Sigma$ is geometrically finite, the horocyclic orbits of $\Omega_h$ are dense or periodic \cite{Hedlund, Ghys, Dalbo}; $h_\R$ is minimal if and only if $\Sigma$ is compact, and otherwise the minimal sets are the periodic orbits. In general, the orbits of $\Omega_h$ are dense, periodic, closed non-compact, or \emph{irregular} (neither dense nor closed). By \cite{DalboStarkov, Bellis-these}, $h_\R(u)$ is dense if and only if $\tilde u(\infty)\in\Lambda_h$, and irregular if and only if $\tilde u(\infty)\in\Lambda_{irr}$. When $\Sigma$ is geometrically infinite, $\Omega_h$ always contains an orbit which is neither compact nor dense \cite{Starkov, DalboStarkov}; $h_\R$ is not minimal, the closed orbits are not minimal, and the only other possible minimal sets are closures of irregular orbits. The closures of irregular orbits understood so far are not minimal \cite{GayeLo, Kulikov, Matsumoto}.

\subsection{Horocyclic convergence.}
We refer to \cite{Beardon, Dalbo} for hyperbolic geometry, the geodesic and horocyclic flows and Busemann cocycles. For $u\in T^1\Sigma$, we denote by $\tilde u$ a lift in $T^1\HH$, $(\tilde u(t))_{t\ge0}$ the geodesic ray it determines, $\tilde u(0)$ its base point and $\tilde u(\infty)\in\partial\HH$ its endpoint. $B_\xi(z,w)$ denotes the Busemann cocycle centred at $\xi$; for $\xi=\infty$,
\[
B_\infty(z,w)=\ln\frac{\Im w}{\Im z}.
\]

\begin{proposition}[{\cite[Ch.~2, Prop.~2.3.2]{Bellis-these}}]\label{prop:bellis}
Let $u,v\in T^1\Sigma$. Then $v\in\cl{h_\R(u)}$ if and only if, for all lifts $\tilde u,\tilde v$, there exists a sequence $(\gamma_n)$ of $\Gamma$ such that (i) $\gamma_n\tilde u(\infty)\to\tilde v(\infty)$ and (ii) $B_{\tilde u(\infty)}(\gamma_n^{-1}i,\tilde u(0))\to B_{\tilde v(\infty)}(i,\tilde v(0))$.
\end{proposition}

Taking $v=g_tu$, $\tilde v=\tilde g_t\tilde u$, so that $\tilde v(\infty)=\tilde u(\infty)$, $\tilde v(0)=\tilde u(t)$, and $\tilde u(0)=i$, we obtain:

\begin{corollary}\label{cor:Tu}
Let $u\in T^1\Sigma$ admit a lift $\tilde u$ based at $i$, and $t\neq0$. Then $g_tu\in\cl{h_\R(u)}$ if and only if there exists a sequence $(\gamma_n)$ of $\Gamma$ such that
\[
\text{(i) } \gamma_n\tilde u(\infty)\to\tilde u(\infty),\qquad \text{(ii) } B_{\tilde u(\infty)}(\gamma_n^{-1}i,i)\to t.
\]
\end{corollary}

\begin{lemma}\label{lem:add}
For every $u\in T^1\Sigma$, the set $T_u$ is stable under addition: $s,t\in T_u$, $s+t\neq0\Rightarrow s+t\in T_u$. In particular, if $t\in T_u$ then $kt\in T_u$ for every $k\ge1$.
\end{lemma}

\begin{proof}
Since $g_t\circ h_s=h_{se^{-t}}\circ g_t$, the map $g_t$ sends $h_\R$-orbits onto $h_\R$-orbits, so $g_t(\cl{h_\R(v)})=\cl{h_\R(g_tv)}$ for every $v$. If $t\in T_u$, then $g_tu\in\cl{h_\R(u)}$, a closed $h_\R$-invariant set, hence $\cl{h_\R(g_tu)}\subset\cl{h_\R(u)}$. If moreover $s\in T_u$, then $g_{s+t}u=g_t(g_su)\in g_t(\cl{h_\R(u)})=\cl{h_\R(g_tu)}\subset\cl{h_\R(u)}$.
\end{proof}

\subsection{Two metric facts.}
The following classical lemma is fundamental in determining the asymptotic thinness.

\begin{lemma}\label{lem:metric}
Let $z=x+iy\in\HH$ and $w\in\HH$.
\begin{enumerate}[label=(\arabic*)]
\item $d(z,w)\ge\left|\ln\dfrac{\Im z}{\Im w}\right|$.
\item If $(0\infty)$ is the vertical geodesic and $z\notin(0\infty)$, then
\[
d((0\infty),z)=2\Argsh\biggl(\sqrt{\frac{|z|-y}{2y}}\biggr)\ge\ln\Bigl(\frac{|z|}{y}-1\Bigr)-\ln2.
\]
\item If $(uv)$ is the geodesic with finite endpoints $u\neq v$ and $z\notin(uv)$, then
\[
d((uv),z)=2\Argsh\biggl(\sqrt{\frac{|z-u||z-v|}{2y|v-u|}-\frac12}\biggr)\ge\ln\Bigl(\frac{|z-u||z-v|}{y|v-u|}-1\Bigr)-\ln2.
\]
\end{enumerate}
\end{lemma}

\begin{proof}
(1) The hyperbolic metric is $ds=\frac{|dz|}{\Im z}\ge\frac{|d(\Im z)|}{\Im z}$. Integrating along the geodesic segment $[z,w]$, we obtain $d(z,w)\ge\int|d\ln\Im z|\ge|\ln\Im z-\ln\Im w|$.

(2) For $z\notin(0\infty)$, we have:
\[
d(z,(0\infty))=\inf_{t\in\R}d(z,ie^t)\quad\text{where}\]\[d(z,ie^t)=2\Argsh\biggl(\frac{|z-ie^t|}{2\sqrt{ye^t}}\biggr)=2\Argsh\biggl(\sqrt{\frac{|z|^2e^{-t}-2y+e^t}{4y}}\biggr).
\]
For $t$ in $\R$, set $f(t)=\dfrac{|z|^2e^{-t}-2y+e^t}{4y}$. The study of $f$ shows that the minimum is attained for $e^t=|z|$ and equals $\dfrac{|z|-y}{2y}$. Hence
\[
d((0\infty),z)=2\Argsh\biggl(\sqrt{\frac{|z|-y}{2y}}\biggr)\ge2\ln\biggl(\sqrt{\frac{|z|-y}{2y}}\biggr)=\ln\Bigl(\frac{|z|}{y}-1\Bigr)-\ln2.
\]

(3) Let $(uv)$ be a geodesic with endpoints $u$ and $v$, and $z\notin(uv)$.

The map $z\mapsto\dfrac{z-u}{(u-v)(z-v)}$ is an isometry of $\HH$ which sends the geodesic $(uv)$ onto the geodesic $(0\infty)$. By point (2),
\[
d((uv),z)=d((0\infty),h(z))=2\Argsh\biggl(\sqrt{\frac{|h(z)|-Im(h(z))}{2Im(h(z))}}\biggr)\ge\ln\Bigl(\frac{|h(z)|}{Im(h(z))}-1\Bigr)-\ln2.
\]
Now $\dfrac{|h(z)|}{Im(h(z))}=\dfrac{|z-u|}{|v-u||z-v|}\dfrac{|z-v|^2}{Im(z)}=\dfrac{|z-u||z-v|}{|u-v|y}$; whence the result.
\end{proof}

\section{Proof of Theorem \ref{thm:main}}

\subsection{Construction of the surface $\Sigma_\delta$.}
Throughout the construction we fix a real number $\delta$ strictly between $1$ and $3$. For every integer $n\ge1$ we set
\[
a_n=\delta+\frac{2n}{n^2+1},\quad b_n=n+(n^2+1)\delta,\quad c_n=\frac1n,\quad d_n=\frac{n^2+1}{n}.
\]
We have $a_nd_n-b_nc_n=1$ and $a_n+d_n>2$, so $h_n(z)=\dfrac{a_nz+b_n}{c_nz+d_n}$ is a hyperbolic isometry of $\HH$.

Let $S_n=\{z\in\HH:\ |(h_n^{-1})'(z)|=1\}$ and $S_{-n}=\{z\in\HH:\ |h_n'(z)|=1\}$ be the isometric circles of $h_n^{-1}$ and of $h_n$. They are Euclidean half-circles of the same radius
\[
R_n=\frac1{c_n}=n
\]
and with centres
\[
X_n=\frac{a_n}{c_n}=na_n=n\delta+\frac{2n^2}{n^2+1}=n\delta+2-\frac{2}{n^2+1},\qquad Y_n=-\frac{d_n}{c_n}=-(n^2+1).
\]
Let $D_n,D_{-n}$ denote the open half-discs bounded by $S_n,S_{-n}$, and $E_{\pm n}=\HH\setminus D_{\pm n}$. Then $h_n(S_{-n})=S_n$ and $h_n(E_{-n})=\cl{D_n}\cap\HH$. Moreover
\[
X_n-n=n(\delta-1)+\frac{2n^2}{n^2+1}>0,\qquad Y_k+k=-(k^2-k+1)<0,
\]
so $\cl{D_n}\subset\{\Re z>0\}$ and $\cl{D_{-k}}\subset\{\Re z<0\}$ (closures in $\cl{\HH}$): the $D_n$ and the $D_{-k}$ have disjoint closures, and none of them meets the imaginary axis. Note also
\[
X_{n+1}-X_n=\delta+\frac{2}{n^2+1}-\frac{2}{(n+1)^2+1}<2n+1,
\]
for $n$ large enough. The idea is to extract an infinite subfamily of half-discs $D_n$ so that $\cl{D_n}\cap\cl{D_p}=\emptyset$ for $n\neq p$. To this end we define the sequence of natural numbers $(p_n)_{n\ge1}$ by the recurrence relation:
\[
p_1=1+E\Bigl(\frac{\delta-1}{2}\Bigr)\ \text{ and }\ p_{n+1}=1+E\Bigl(\frac{(\delta+1)p_n}{\delta-1}\Bigr)=1+p_n+E\Bigl(\frac{2p_n}{\delta-1}\Bigr)\ \ \forall\, n\ge1,
\]
where $E$ denotes the integer part. Thus $(p_n)_{n\ge1}$ is a strictly increasing sequence of integers and satisfies:
\[
p_{n+1}-p_n=1+E\Bigl(\frac{2p_n}{\delta-1}\Bigr)\ge1+E\Bigl(\frac{2p_1}{\delta-1}\Bigr)\ge1+E\Bigl(\frac{2}{\delta-1}+\frac{2}{\delta-1}E\Bigl(\frac{\delta-1}{2}\Bigr)\Bigr)=2.
\]
\[
p_{n+1}=1+p_n+E\Bigl(\frac{2p_n}{\delta-1}\Bigr)>1+p_n+\frac{2p_n}{\delta-1}-1>\Bigl(1+\frac{2}{\delta-1}\Bigr)p_n>2p_n\ \text{ since } 1<\delta<3.
\]

\begin{lemma}\label{lem:disjoint}
The closed half-discs $\cl{D_{p_n}}$, $\cl{D_{-p_n}}$ $(n\ge1)$ are pairwise disjoint in $\cl{\HH}$.
\end{lemma}

\begin{proof}
On the one hand, for every $n\ge1$ we have:
\[
X_{p_{n+1}}-X_{p_n}-p_{n+1}-p_n=\frac{2(p_{n+1}^2-p_n^2)}{(p_{n+1}^2+1)(p_n^2+1)}+(\delta-1)p_{n+1}-(\delta+1)p_n.
\]
Now,
\[
(\delta-1)p_{n+1}=(\delta-1)\biggl(1+E\Bigl(\frac{(\delta+1)p_n}{\delta-1}\Bigr)\biggr)>(\delta-1)\Bigl(\frac{(\delta+1)p_n}{\delta-1}\Bigr)=(\delta+1)p_n.
\]
Consequently we obtain:
\[
X_{p_{n+1}}-X_{p_n}-p_{n+1}-p_n\ge(\delta-1)p_{n+1}-(\delta+1)p_n>(\delta+1)p_n-(\delta+1)p_n=0
\]
and therefore the half-discs $\cl{D}_{p_n}$ and $\cl{D}_{p_{n+1}}$ are disjoint.

On the other hand, for every $n\ge1$, we also have:
\[
Y_{p_n}-Y_{p_{n+1}}=p_{n+1}^2-p_n^2=(p_{n+1}+p_n)(p_{n+1}-p_n)>p_{n+1}+p_n
\]
and therefore the half-discs $\cl{D}_{-p_n}$ and $\cl{D}_{-p_{n+1}}$ are also disjoint.
\end{proof}

Let $\Gamma=\langle h_{p_n};\ n\ge1\rangle$ and $\Delta=\bigcap_{n\ge1}E_{p_n}\cap E_{-p_n}$. By Lemma~\ref{lem:disjoint} and the classical Klein--Maskit combination argument (ping-pong) \cite{Beardon}, $\Gamma$ is a free Fuchsian group (an infinitely generated Schottky group) with fundamental domain $\Delta$. Every non-trivial element is hyperbolic: a reduced word involves only finitely many generators, hence belongs to a Schottky group of finite type with pairwise disjoint closed discs. Thus $\Gamma$ has neither elliptic nor parabolic elements. Being infinitely generated, $\Gamma$ is geometrically infinite. We set $\Sigma_\delta=\Gamma\backslash\HH$. Note that the imaginary axis $i\R_+^*$ is contained in the interior of $\Delta$.

\begin{figure}[ht]
\centering
\includegraphics[width=0.82\textwidth]{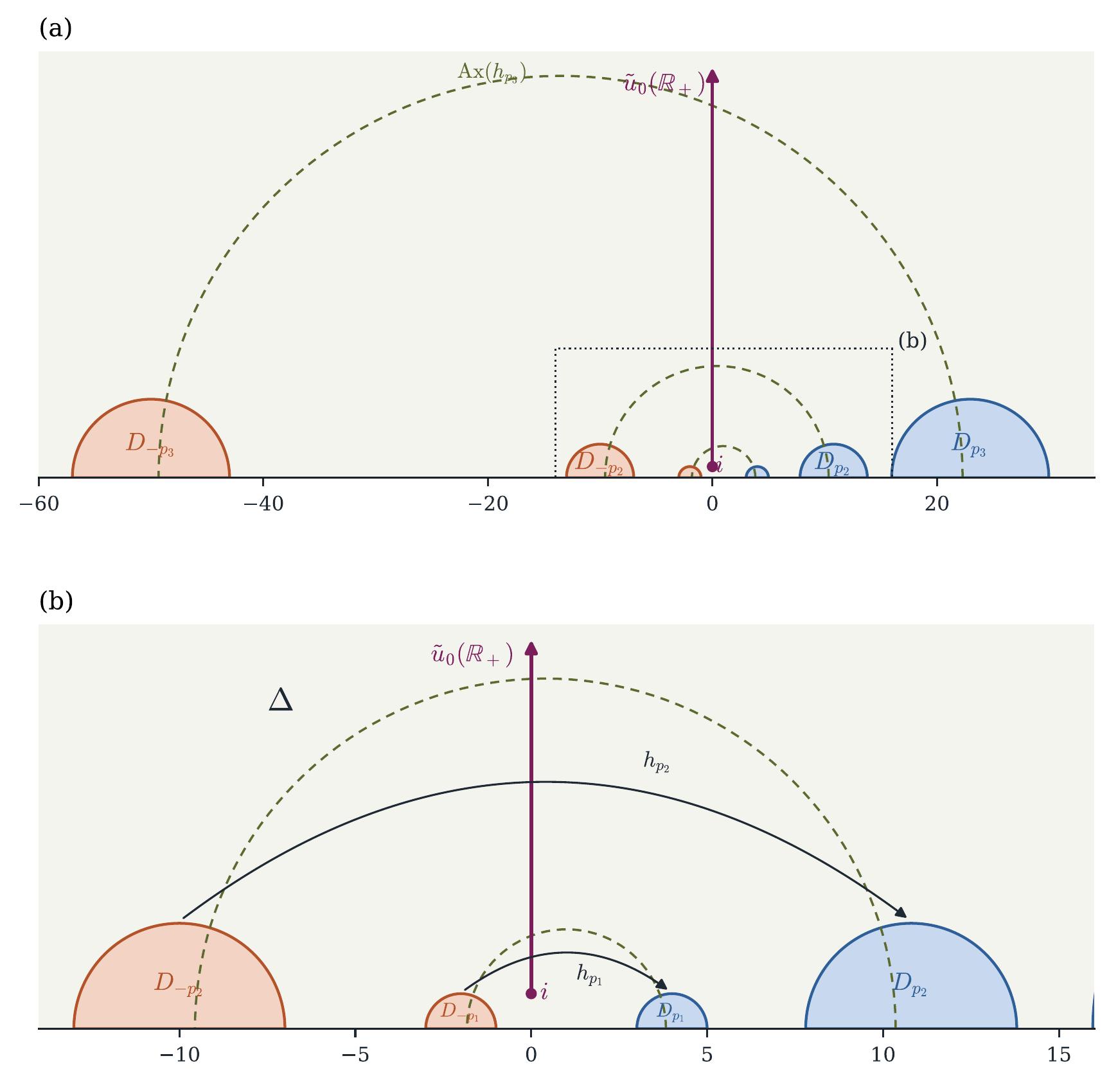}
\caption{The fundamental domain $\Delta$ (light background) is the complement of the half-discs $D_{\pm p_n}$, and $h_{p_n}$ sends $S_{-p_n}$ onto $S_{p_n}$. Dashed: the axes $\mathrm{Ax}(h_{p_n})$. The vertical ray from $i$ is the lift $\tilde u_0(\R_+)$. (a) Overall view; (b) enlargement near the origin.}
\end{figure}

\paragraph{\textbf{Ping-pong.}} Let $\mathcal A=\{h_{p_n}^{\pm1};\ n\ge1\}$. For $s\in\mathcal A$ we define its \emph{target disc} $T(s)=D_p$ if $s=h_p$, $T(s)=D_{-p}$ if $s=h_p^{-1}$. The isometric circle of $s$ bounds $T(s^{-1})$, and $s$ sends $\HH\setminus\cl{T(s^{-1})}$ into $T(s)$. By (3.1), $|s'(z)|<1$ for $z\notin\cl{T(s^{-1})}$.

\begin{lemma}\label{lem:pingpong}
Let $\gamma=s_1s_2\cdots s_m$ be a reduced word ($s_j\in\mathcal A$, $s_{j+1}\neq s_j^{-1}$), $m\ge1$, and let $z\in\Delta^\circ$. Set $w_{m+1}=z$ and $w_j=s_j(w_{j+1})$. Then $w_j\in T(s_j)$ for $1\le j\le m$, $|s_j'(w_{j+1})|<1$, and consequently $|\gamma'(z)|<1$.
\end{lemma}

\begin{proof}
$w_{m+1}=z$ lies outside all the closed discs, so $w_m=s_m(z)\in T(s_m)$ and $|s_m'(z)|<1$. If $w_{j+1}\in T(s_{j+1})$, then $T(s_{j+1})\neq T(s_j^{-1})$ since the word is reduced, so by Lemma~\ref{lem:disjoint}, $w_{j+1}\notin\cl{T(s_j^{-1})}$; hence $w_j=s_j(w_{j+1})\in T(s_j)$ and $|s_j'(w_{j+1})|<1$. We conclude by the chain rule.
\end{proof}

\subsection{The irregular orbit $h_\R(u_0)$.}
Since $h_{p_n}(i)\in D_{p_n}$ and $X_{p_n}-p_n\to+\infty$, we have $h_{p_n}(i)\to\infty$, so $\infty\in\Lambda$. Let $\tilde u_0\in T^1\HH$ be the vector based at $i$ and directed towards $\infty$, and $u_0$ its projection in $T^1\Sigma_\delta$. Then $u_0\in\Omega_h$.

\begin{proposition}\label{prop:irr}
The orbit $h_\R(u_0)$ is irregular, $2\ln\delta\in T_{u_0}$, and $\cl{h_\R(u_0)}$ is not $h_\R$-minimal.
\end{proposition}

\begin{proof}
$\infty\notin\Lambda_p$: $\Gamma$ has no parabolic element.

$\infty\notin\Lambda_d$: Since ${h_{p_n}}^{-1}(z)=\frac{d_{p_n}z-b_{p_n}}{-c_{p_n}z+a_{p_n}}$, (3.1) gives
\[
\Im h_{p_n}^{-1}(i)=\frac{1}{c_{p_n}^2+a_{p_n}^2}=\delta^2+\frac{4\delta p_n}{p_n^2+1}+\frac{4p_n^2}{(p_n^2+1)^2}+\frac{1}{p_n^2}\xrightarrow[p_n\to\infty]{}\delta^2.
\]
The elements $h_{p_n}^{-1}$ are pairwise distinct, so $\Im h_{p_n}^{-1}(i)\to\delta^{-2}>0$ and Remark~\ref{rem:horo} applies.

$\infty\notin\Lambda_h$: let $\gamma\neq\Id$. Since $i\in\mathring{\Delta}$, Lemma~\ref{lem:pingpong} and (3.1) give $\Im\gamma(i)=|\gamma'(i)|<1$. Thus $\Im\gamma(i)\le1$ for every $\gamma\in\Gamma$, and Remark~\ref{rem:horo} applies.

Hence $\infty\in\Lambda_{irr}$ and $h_\R(u_0)$ is irregular. Now take $\gamma_n=h_{p_n}$: $\gamma_n(\infty)=X_{p_n}\to\infty$ and
\[
B_\infty(\gamma_n^{-1}i,i)=-\ln\Im h_{p_n}^{-1}(i)=\ln(a_{p_n}^2+c_{p_n}^2)\longrightarrow\ln\delta^2=2\ln\delta\neq0.
\]
By Corollary~\ref{cor:Tu}, $2\ln\delta\in T_{u_0}$. Since $h_\R(u_0)$ is not closed and $T_{u_0}\neq\emptyset$, $\cl{h_\R(u_0)}$ is not minimal \cite{Bellis, Bellis-these}.
\end{proof}

\subsection{Asymptotic thinness of $u_0(\R_+)$.}

\begin{definition}\label{def:inj}
For $x\in\Sigma_\delta$ with lift $\tilde x$, $\Inj(x)=\inf_{\gamma\in\Gamma\setminus\{\Id\}}d(\tilde x,\gamma\tilde x)$. For $u\in T^1\Sigma_\delta$, the asymptotic thinness of $u(\R_+)$ is $\Inj(u(\R_+))=\liminf_{t\to+\infty}\Inj(u(t))$.
\end{definition}

The lift of $u_0(\R_+)$ is $z_t=ie^t$, $t\ge0$. We first estimate the distance from $z_t$ to the circles $S_{-p_n}$.

\begin{lemma}\label{lem:dist}
There exists a real number $T>0$ such that for all $t\ge T$ and every integer $p_n\ge1$,
\[
d(ie^t,S_{-p_n})\ge\frac t2-\ln3.
\]
\end{lemma}

\begin{proof}
\emph{Step 1: $S_{-p_n}$ is the isometric circle $\{|c_{p_n}z+d_{p_n}|=1\}$ of $h_{p_n}$.}
\[
|c_{p_n}z+d_{p_n}|=1\iff\Bigl|z+\frac{d_{p_n}}{c_{p_n}}\Bigr|=\frac1{c_{p_n}}\iff|z-Y_{p_n}|=p_n,\qquad Y_{p_n}=-(p_n^2+1).
\]
Thus $S_{-p_n}$ is the half-circle with centre $Y_{p_n}$ and radius $p_n$, with endpoints:
\[
u_n=Y_{p_n}+p_n=-(p_n^2-p_n+1)\ \text{ and }\ v_n=Y_{p_n}-p_n=-(p_n^2+p_n+1).
\]
Consequently:
\[
u_nv_n=(p_n^2-p_n+1)(p_n^2+p_n+1)=(p_n^2+1)^2-p_n^2=p_n^4+p_n^2+1\ \text{ and }\ |u_n-v_n|=2p_n.
\]

\emph{Step 2: the distance.} Set $z=ie^t$. By Lemma~\ref{lem:metric}, we have:
\[
d(ie^t,S_{p_n})\ge\ln\biggl(\frac{|ie^t-u_n||ie^t-v_n|}{e^t|v_n-u_n|}-1\biggr)-\ln2.
\]
\[
\frac{|ie^t-u_n||ie^t-v_n|}{e^t|v_n-u_n|}=\frac{\sqrt{(-e^{2t}+u_nv_n)^2+e^{2t}(v_n+u_n)^2}}{2p_ne^t}=\frac{\sqrt{(-e^t+u_nv_ne^{-t})^2+(v_n+u_n)^2}}{2p_n}.
\]
Now,
\[
(-e^t+u_nv_ne^{-t})^2+(v_n+u_n)^2=e^{2t}+u_n^2v_n^2e^{-2t}+u_n^2+v_n^2\ge(e^t+u_nv_ne^{-t})^2.
\]
We then obtain:
\[
\frac{|ie^t-u_n||ie^t-v_n|}{e^t|v_n-u_n|}\ge\frac{e^t+u_nv_ne^{-t}}{2p_n}=\frac{e^t+(p_n^4+p_n^2+1)e^{-t}}{2p_n}
\]
and therefore
\[
d(ie^t,S_{-p_n})\ge\ln\biggl(\frac{e^t+(p_n^4+p_n^2+1)e^{-t}}{2p_n}-1\biggr)-\ln2.
\]

\emph{Step 3: A lower bound uniform in $p_n$.} Dropping the term $p_n^2+1>0$.

Set $\varphi(x)=:\frac{e^t}{2x}+\frac{x^3}{2e^t}$, $x>0$. We have $\varphi'(x)=-\frac{e^t}{2x^2}+\frac{3x^2}{2e^t}$, which vanishes only at $x_*=3^{-1/4}e^{t/2}$. Since $\varphi(x)\to+\infty$ as $x\to0^+$ and as $x\to+\infty$, this unique critical point is the global minimum on $]0,+\infty[$, namely
\[
\varphi(x_*)=\frac{3^{1/4}}{2}e^{t/2}+\frac{3^{-3/4}}{2}e^{t/2}=\frac{e^{t/2}}{2}\,3^{-3/4}(3+1)=\frac{2}{3^{3/4}}\,e^{t/2}.
\]
The lower bound therefore holds for every real $x>0$, in particular for every integer $p_n\ge1$.

\emph{Consequence.} Since $t>\ln\Bigl(\dfrac{\sqrt{27}}{4}\Bigr)$ implies that $\dfrac{2e^{t/2}}{3^{3/4}}>1$, then
\[
d(ie^t,S_{-p_n})\ge\ln\Bigl(\frac{2e^{t/2}}{3^{3/4}}-1\Bigr)-\ln2=\frac t2-\frac{3\ln3}{4}+\ln\Bigl(1-\frac{3^{3/4}e^{-t/2}}{2}\Bigr).
\]
Since $\ln\Bigl(1-\dfrac{3^{3/4}e^{-t/2}}{2}\Bigr)$ tends to $0$ as $t$ tends to $+\infty$, there exists $t_0>0$ such that for $t>t_0$ we have $\ln\Bigl(1-\dfrac{3^{3/4}e^{-t/2}}{2}\Bigr)\ge\dfrac{-\ln3}{4}$. Setting $T=\max\Bigl\{t_0,\ln\Bigl(\dfrac{\sqrt{27}}{4}\Bigr)\Bigr\}$ we then obtain:
\[
\forall\, t\ge T,\ d(ie^t,S_{-p_n})\ge\frac t2-\ln3.
\]
\end{proof}

\begin{proposition}\label{prop:inj}
There exists a real number $A>0$ such that for all $t\ge A$, $\Inj(u_0(t))\ge\dfrac t2-\ln3$; in particular $\Inj(u_0(\R_+))=+\infty$.
\end{proposition}

\begin{proof}
\emph{Set-up.} The lift of $u_0(\R_+)$ is $z_t=ie^t$, parametrised by arc length since $d(i,ie^t)=t$. So $\Inj(u_0(t))=\inf_{\gamma\neq\Id}d(z_t,\gamma z_t)$. Write $y=e^t$ and $z=z_t$. Since $\Re z=0$ and each closed half-disc $\cl{D_{\pm p_n}}$ is contained in $\{\Re>0\}$ or in $\{\Re<0\}$, we have $z\in\mathring{\Delta}$. Let $\gamma=s_1\cdots s_m\neq\Id$ be a reduced word, $s_j\in\{h_{p_n}^{\pm1},n\ge1\}$. By Lemma~\ref{lem:pingpong}, $\gamma(z)\in T(s_1)$. Three cases cover all $\gamma\neq\Id$.

\emph{Case A: $s_1=h_{p_n}^{-1}$.} Then $\gamma z\in D_{-p_n}$ whereas $z\notin\cl{D_{-p_n}}$. The geodesic segment $[z,\gamma z]$ is connected, starts outside and ends inside $\cl{D_{-p_n}}$; it therefore meets $S_{-p_n}$ at a point $\zeta$. By Lemma~\ref{lem:dist}, for $t>T$, we have:
\[
d(z,\gamma z)\ge d(z,\zeta)\ge d(z,S_{-p_n})>\frac t2-\ln3.
\]

\emph{Case B: $\gamma=h_{p_n}$.}
We have $\Im h_{p_n}(iy)=\frac{y}{|c_{p_n}iy+d_{p_n}|^2}=\frac{y}{c_{p_n}^2y^2+d_{p_n}^2}$. Since $c_{p_n}^2y^2+d_{p_n}^2\ge2yc_{p_n}d_{p_n}$, then:
\[
\frac{\Im z}{\Im h_{p_n}(z)}=c_{p_n}^2y^2+d_{p_n}^2\ge2c_{p_n}d_{p_n}y=2\,\frac{{p_n}^2+1}{p_n^2}\,y\ge2y,
\]
and point (1) of Lemma~\ref{lem:metric} gives $d(z,h_{p_n}z)\ge\ln(2y)=t+\ln2$.

\emph{Case C: $s_1=h_{p_n}$ and $m\ge2$.}
Set $w=s_2\cdots s_m(z)$. Lemma~\ref{lem:pingpong} applied to the reduced word $s_2\cdots s_m$ gives $w\in T(s_2)$. Since the word is reduced, $s_2\neq h_{p_n}^{-1}$, so $T(s_2)\neq D_{-p_n}$. Moreover $\gamma z=h_{p_n}(w)$. Since $c_{p_n}w+d_{p_n}=c_{p_n}(w-Y_{p_n})$ and $c_{p_n}=1/p_n$, we have:
\begin{equation}\tag{3.4}\label{eq:34}
\Im\gamma z=\Im(h_{p_n}(w))=\frac{\Im w}{|c_{p_n}w+d_{p_n}|^2}=\frac{p_n^2\,\Im w}{|w-Y_{p_n}|^2}.
\end{equation}
\begin{enumerate}[label=($c_{\arabic*}$),leftmargin=2.2em]
\item $T(s_2)=D_{q_k}$, $q_k$ arbitrary (possibly $q_k=p_n$). Since $w\in D_{q_k}$, $\Im w<q_k$ and
\[
\Re w>X_{q_k}-q_k=q_k(\delta-1)+\frac{2q_k^2}{q_k^2+1}>q_k(\delta-1).
\]
Hence $|w-Y_{p_n}|^2\ge(\Re w-Y_{p_n})^2>(q_k(\delta-1)+p_n^2+1)^2>(q_k(\delta-1)+p_n^2)^2\ge2q_kp_n^2(\delta-1)$. We obtain:
\[
\Im\gamma z<\frac{p_n^2q_k}{4q_k(\delta-1)p_n^2}=\frac{1}{4(\delta-1)}.
\]
\item $T(s_2)=D_{-q_k}$ with $q_k\neq p_n$. Since $w\in D_{-q_k}$, $|w-Y_{q_k}|<q_k$ and $\Im w<q_k$, and by the reverse triangle inequality
\[
|w-Y_{p_n}|\ge|Y_{q_k}-Y_{p_n}|-|w-Y_{q_k}|>|q_k^2-p_n^2|-q_k.
\]
Recall that for every $j\ge1$, $p_{j+1}\ge2p_j+1$.

- If $q_k>p_n$, then $q_k\ge p_{n+1}\ge2p_n+1\ge3$ and $p_n\le\frac{q_k-1}{2}<\dfrac{q_k}{2}$, so
\[
q_k^2-p_n^2-q_k\ge q_k^2-\frac{(q_k-1)^2}{4}-q_k=\frac{(3q_k+1)(q_k-1)}{4}\ge\frac{q_k^2}{2},
\]
since $q_k-1\ge\frac{2q_k}{3}$ and $3q_k+1>3q_k$ and $q_k\ge3$. By \eqref{eq:34},
\[
\Im\gamma z<\frac{p_n^2q_k}{(q_k^2/2)^2}=\frac{4p_n^2}{q_k^3}\le\frac{4(q_k/2)^2}{q_k^3}=\frac{1}{q_k}\le\frac13.
\]

- If $q_k<p_n$, then $p_n\ge2q_k+1\ge3$ and $q_k\le\frac{p_n-1}{2}$, so
\[
p_n^2-q_k^2-q_k\ge p_n^2-\frac{(p_n-1)^2}{4}-\frac{p_n-1}{2}=\frac{3p_n^2+1}{4}\ge\frac34p_n^2.
\]
By \eqref{eq:34},
\[
\Im\gamma z<\frac{p_n^2q_k}{(3p_n^2/4)^2}=\frac{16\,q_k}{9p_n^2}\le\frac{16\,(p_n/2)}{9p_n^2}=\frac{8}{9p_n}\le\frac{8}{27}.
\]
\end{enumerate}
In all sub-cases $\Im\gamma z\le\max\Bigl(\dfrac{1}{4(\delta-1)},\dfrac13\Bigr)$. Setting $K=\min(3,4(\delta-1))$, point (1) of Lemma~\ref{lem:metric} gives
\[
d(z,\gamma z)\ge t+\ln K=\frac t2-\ln3+\frac t2+\ln3+\ln K>\frac t2-\ln3
\]

\emph{Conclusion.} For every $\gamma\neq\Id$, for $t>T$, $d(z_t,\gamma z_t)\ge\dfrac t2-\ln3$ and therefore $\Inj(u_0(t))\ge\dfrac t2-\ln3$ for every $t\ge T$. Hence $\Inj(u_0(\R_+))=\liminf_{t\to+\infty}\Inj(u_0(t))=+\infty$.
\end{proof}

\emph{Why does case C require a subtler argument?} For a word of length at least $2$ with $s_1=h_{p_n}$, the distance from $z$ to $S_{p_n}$ is not enough: by point (3) of Lemma~\ref{lem:metric}, $d(ie^t,S_{p_n})=\Argsh\Biggl(\sqrt{\dfrac{\bigl(e^t+(X_{p_n}^2-p_n^2)e^{-t}\bigr)^2}{2p_n}-1}\Biggr)$, which for $e^t=p_n$ tends to $\Argsh\bigl(\frac{\delta^2-1}{2}\bigr)$ and remains bounded. This is why we bound the height $\Im\gamma z$ instead. The fact that $p_{n+1}>2p_n$ is indispensable only in sub-case $(c_2)$; otherwise (for example with $q_k=p_{n+1}$), we would only have $|w-Y_{p_n}|\gtrsim p_n$, and $\Im\gamma z$ could be of the order of $p_n\approx e^t$.

\subsection{The set $T_{u_0}$.}

\begin{proposition}\label{prop:Tu}
$T_{u_0}$ contains the sequence $t_k=2k\ln\delta$, $k\ge1$, which tends to $+\infty$.
\end{proposition}

\begin{proof}
By Proposition~\ref{prop:irr}, $2\ln\delta\in T_{u_0}$, with $2\ln\delta>0$ since $\delta>1$. Lemma~\ref{lem:add} gives $2k\ln\delta\in T_{u_0}$ for every $k\ge1$.
\end{proof}

This completes the proof of Theorem~\ref{thm:main}.\hfill$\square$


\begin{thebibliography}{11}
\bibitem{Beardon} A.~F. Beardon, \emph{The geometry of discrete groups}, Springer-Verlag, New York, 1983.
\bibitem{Bellis-these} A.~Bellis, \emph{Étude topologique du flot horocyclique : le cas des surfaces géométriquement infinies}, PhD thesis, Université de Rennes~1.
\bibitem{Bellis} A.~Bellis, On the links between horocyclic and geodesic orbits on geometrically infinite surfaces, \emph{J. Éc. polytech. Math.} 5 (2018), 443--454.
\bibitem{Dalbo} F.~Dal'Bo, \emph{Trajectoires géodésiques et horocycliques}, Savoirs Actuels, EDP Sciences--CNRS, 2007.
\bibitem{DalboStarkov} F.~Dal'Bo, A.~N. Starkov, On a classification of limit points of infinitely generated Schottky groups, \emph{J. Dyn. Control Syst.} 6 (2000), no.~4, 561--578.
\bibitem{GayeLo} M.~Gaye, C.~Lo, Sur l'inexistence d'ensembles minimaux pour le flot horocyclique, \emph{Confluentes Math.} 9 (2017), no.~1, 95--104.
\bibitem{Ghys} É.~Ghys, Dynamique des flots unipotents sur les espaces homogènes, Séminaire Bourbaki 1991--92, \emph{Astérisque} 206 (1992), 93--136.
\bibitem{Hedlund} G.~A. Hedlund, Fuchsian groups and transitive horocycles, \emph{Duke Math. J.} 2 (1936), no.~3, 530--542.
\bibitem{Kulikov} M.~Kulikov, The horocycle flow without minimal sets, \emph{C. R. Acad. Sci. Paris, Sér.~I} 338 (2004), 477--480.
\bibitem{Matsumoto} S.~Matsumoto, Horocycle flows without minimal sets, \emph{J. Math. Sci. Univ. Tokyo} 23 (2016), no.~3, 661--673.
\bibitem{Starkov} A.~N. Starkov, Fuchsian groups from the dynamical viewpoint, \emph{J. Dyn. Control Syst.} 1 (1995), no.~3, 427--445.
\end{thebibliography}
\end{document}